\documentclass[11pt]{article}
\usepackage{amsmath,amssymb,amsthm}
\usepackage[margin=1in]{geometry}
\usepackage{hyperref}
\usepackage{booktabs}

\newtheorem{theorem}{Theorem}
\newtheorem{proposition}[theorem]{Proposition}
\newtheorem{lemma}[theorem]{Lemma}

\newtheorem{example}[theorem]{Example}
\newtheorem{openproblem}[theorem]{Open Problem}
\newtheorem*{remark}{Remark}

\begin{document}

\title{The Family $X_{qp^s}$ for Distinct Primes $q,p$ in
Connection with Carmichael's Totient Conjecture}

\author{Tian-Xiao He 
\\
{\small Department of Mathematics}\\
 {\small Illinois Wesleyan University}\\
 {\small Bloomington, IL 61702-2900, USA}\\
}

\date{}

%%%%%%%%%%%%%%%%%%%%%%%%%%%%%%%%%%%%%%%%%%%%%%%%%%%%%%%%%%%%%%%%%%%%%%
\maketitle
\setcounter{page}{1}
\pagestyle{myheadings}
\markboth{T. X. He }
{Family $X_{qp^s}$ in
Connection with Carmichael's Totient Conjecture}

%%%%%%%%%%%%%%%%%%%%%%%%%%%%%%%%%%%%%%%%%%%%%%%%%%%%%%%%%%%%%%%%%%%%%%
\begin{abstract}
\noindent 
Shannon, Shiue, Saito and the author \cite{SSHS2025} studied the sets $X_k=\{n\in\mathbb{N}:\varphi(n)=k\}$ in
connection with Carmichael's totient conjecture, and proved that for an odd prime $p\neq3$ and
$s\in\{1,2,3,4\}$, $X_{2p^s}$ is empty whenever $2p^s+1$ is not prime (with $s=1$ fully characterized in both
directions), conjecturing that this persists for every $s\ge1$. We resolve that conjecture in full, obtaining
a complete characterization of $X_{2p^s}$ for every $s\ge1$ and every odd prime $p\neq3$, together with a
parity dichotomy: $X_{2p^s}=\varnothing$ unconditionally whenever $s$ is even, while for odd $s$ both outcomes
genuinely occur, determined exactly by the primality of $2p^s+1$. We then close the one remaining gap at
$p=3$, where the proof technique degenerates, with a separate, unconditional characterization of $X_{2\cdot3^s}$.

Building on these results, we consider the natural generalization $X_{qp^s}$, where $q\neq p$ are distinct
primes and $s\ge1$. We show that this generalization splits into three cases of very different character. If
$q$ and $p$ are both odd, the problem is trivial: $qp^s$ is odd, so $X_{qp^s}=\varnothing$ by elementary
parity, with no case-analysis needed. If $q=2$, the problem is exactly the family $X_{2p^s}$ resolved above.
If $p=2$, giving the ``mirror'' family $X_{2^sq}$ for an odd prime $q$, the problem is genuinely different.
Beyond $s=1$, we give a complete resolution for $s=2$, recovering, with an explicit proof and solution sets,
the qualitative dichotomy of \cite[Thm.~2.3.1]{SSHS2025}, but the problem remains open for $s\ge3$. We explain
precisely why the proof technique that resolves $X_{2p^s}$ does not transfer, and we report computational
evidence showing that $|X_{2^sq}|$ does \emph{not} obey the clean ``$0$ or $2$'' dichotomy of the other two
cases once $s\ge3$, even under the congruence condition that regularizes $s=2$.


\end{abstract}

\section{Introduction}
\label{sec:intro}

Let $\varphi$ denote Euler's totient function, and for a positive integer $k$ write
$$
X_k = \{\, n \in \mathbb{N} : \varphi(n) = k \,\}.
$$
In 1907, Carmichael claimed to have \emph{proved} that $|X_k|\neq1$ for every $k$ in the image of $\varphi$ ---
that is, that no value of $\varphi$ is ever attained by a unique $n$ \cite{Carmichael1907}. He tabulated
solutions of $\varphi(m)=k$ for $k\le1000$ the following year \cite{Carmichael1908}, and the statement appears
again, as an exercise, in his 1914 textbook \cite[Ch.~2, Ex.~8]{Carmichael1914}. The original 1907 argument,
however, was flawed; Carmichael himself retracted the claim in 1922, restating it as a \emph{conjecture} and
showing that any counterexample must exceed $10^{37}$ \cite{Carmichael1922}. It has remained open ever since,
and is now universally known as \emph{Carmichael's (totient) conjecture}; it appears as Problem B39 in Guy's
\emph{Unsolved Problems in Number Theory}, alongside the closely related Problem B36 \cite[B36, B39]{Guy2004}.

The distribution of $|X_k|$ as $k$ ranges over totient values was studied in depth by Ford, who determined the
typical and extremal order of the number of solutions of $\varphi(n)=m$ \cite{Ford1998,Ford1999}; this body of
work underlies the modern analytic approach to Carmichael's conjecture and to the two closely related questions
below.

Erd\"os sharpened the uniqueness question in a different direction, asking (Problem \#51 of the Erd\"os Problems
database \cite{Erdos51}) whether there exists an infinite set $A\subset\mathbb{N}$ of totient values such that,
writing $n_a$ for the \emph{smallest} $n$ with $\varphi(n)=a$, one has $n_a/a\to\infty$ as $a\to\infty$ through
$A$; the same page records that Erd\"os proved that if some $t$ has $|X_t|=1$, then in fact infinitely many $t$
do, tracing the question directly back to Carmichael's original uniqueness claim. A closely related extremal
question, Problem \#694 \cite{Erdos694}, asks for $\max_{n\le x} f_{\max}(n)/f_{\min}(n)$, where $f_{\max}(n)$
and $f_{\min}(n)$ denote the largest and smallest preimages of a totient value $n$. As this note was being
prepared, that page was updated to record a proof that
$$
\max_{n\le x} \frac{f_{\max}(n)}{f_{\min}(n)} = (e^{\gamma}+o(1))\log\log x,
$$
attributed on the site to the AI system GPT-5.5 Pro (prompted by L.~Price) \cite{Erdos694}. We note explicitly,
since the question is adjacent to ours, that this is at present a crowd-maintained website claim rather than a
peer-reviewed publication; and --- more importantly for the present note --- that its proof is of a
fundamentally different, asymptotic and analytic, character from the elementary case analysis carried out
below. That result concerns the \emph{typical and extremal} behavior of preimage sizes as $n$ ranges over
\emph{all} totient values up to $x$, built on Ford's distributional machinery \cite{Ford1998,Ford1999}, whereas
the results of this note are \emph{exact, unconditional characterizations} of $X_k$ for specific infinite
families of structured $k$ (such as $k=2p^s$ or $k=2^sq$), obtained by elementary divisibility arguments with
no appeal to density or asymptotics. Beyond the shared ancestry in Carmichael's and Erd\"os's original
multiplicity questions, we see no direct technical overlap between the two approaches; we record the
connection here so that the reader is aware we have considered it.

Our own starting point is \cite{SSHS2025}, where it is shown that for an odd prime $p\neq3$ and
$s\in\{1,2,3,4\}$, $X_{2p^s}=\varnothing$ whenever $2p^s+1$ is not prime (with $s=1$ fully characterized in both
directions), and it is conjectured that this persists for every $s\ge1$. In \S3 below we resolve that
conjecture in full, giving a complete characterization of $X_{2p^s}$ for every $s\ge1$ and every odd prime
$p\neq3$, together with the accompanying parity dichotomy. We then contribute two further pieces of this
picture: we close the gap left open by the standing hypothesis $p\neq3$ (\S\ref{sec:p3}), and we investigate
the natural companion family obtained by exchanging the roles of the two primes, $X_{2^sq}$ with $q$ an odd
prime, which we resolve completely for $s\le2$ and show --- with supporting computation --- to require
genuinely new ideas for $s\ge3$ (\S\ref{sec:mirror}).

We record three basic facts about $\varphi$, used throughout, all due to \cite{SSHS2025}:

\begin{itemize}
\item[(P1)] (\cite[Prop.~1.1.1]{SSHS2025}) If $q$ is prime and $a \geq 1$ is an integer, then
$\varphi(q^a) = q^{a-1}(q-1)$. Separately, by direct verification, $\varphi(q^0) = \varphi(1) = 1$ (the only
positive integer $\le 1$ coprime to $1$ is $1$ itself).
\item[(P2)] (\cite[Prop.~1.1.6]{SSHS2025}) If $m \mid n$ then $\varphi(m) \mid \varphi(n)$.
\item[(P3)] (\cite[Prop.~1.1.8]{SSHS2025}) For every positive integer $n$, $\varphi(n)$ is either $1$ or even.
\end{itemize}

Throughout this paper, we fix distinct primes $q\neq p$ and an integer $s\ge1$, set $k=qp^s$, and discuss the construction 
of $X_k$. The remainder of this note is organized as follows. \S2 disposes of the trivial case where $q$ and $p$ are both
odd primes, via (P3). \S3 gives the complete resolution of the case $q=2$: the master formula for $X_{2p^s}$
(Proposition~\ref{prop:master}), valid for every odd prime $p\neq3$ and every $s\ge1$, together with its
even/odd parity dichotomy (Theorems~\ref{thm:even}, \ref{thm:odd}). \S\ref{sec:p3} closes the one remaining gap
in that family, $p=3$, with a complete and
unconditional characterization of $X_{2\cdot3^s}$ (Proposition~\ref{prop:p3-master}). \S\ref{sec:mirror} turns
to the mirror family $X_{2^sq}$ obtained by exchanging the roles of $q$ and $p$: we resolve it completely for
$s\le2$ (Proposition~\ref{prop:mirror-master}), document computationally how the clean dichotomy breaks down
for $s\ge3$, and record the resulting open problem, closing with a remark on the shared boundary $q=p=2$ where
the resolved and mirror families coincide.

\section{The trivial case of $X_{q\cdot p^s}$ with $q$ and $p$ both odd}

\begin{proposition}\label{prop:trivial}
If $q$ and $p$ are distinct odd primes and $s\ge1$, then $X_{qp^s}=\varnothing$.
\end{proposition}

\begin{proof}
Since $q$ and $p$ are odd, $qp^s$ is a product of odd numbers, hence odd. Also $qp^s \ge 3\cdot5=15>1$. By (P3),
$\varphi(n)$ is never an odd number greater than $1$, so no $n$ satisfies $\varphi(n)=qp^s$. Hence
$X_{qp^s}=\varnothing$.
\end{proof}

\begin{remark}
This case requires none of the machinery below: it is a one-line parity obstruction, not a consequence of the
finer structure of $X_k$.
\end{remark}

\section{The case of $X_{2p^s}$ for $q=2$}

This section solves $X_{2p^s}$ for every odd prime $p\neq3$ and every $s\ge1$ (extending Theorems 2.2.1--2.2.4 of
\cite{SSHS2025}, which had treated only $s\le4$).

We give some standing hypotheses for this section. $p$ is an odd prime with $p\neq3$ (i.e.\ $p\ge5$), $s\ge1$, and
$q_j:=2p^j+1$ for $j=1,\dots,s$. 

\begin{lemma}\label{lem:p-not-divide}
No element of $X_{2p^s}$ is divisible by $p$.
\end{lemma}

\begin{proof}
Suppose $n \in X_{2p^s}$ and $p \mid n$. By (P2) applied to $m = p$, $\varphi(p) = p-1$ divides $\varphi(n) = 2p^s$.
Since $p$ and $p-1$ are consecutive integers, $\gcd(p-1,p) = 1$, and hence $\gcd(p-1, p^s) = 1$. Because
$(p-1) \mid 2p^s$ and $\gcd(p-1,p^s)=1$, it follows that $(p-1) \mid 2$. Thus $p - 1 \in \{1,2\}$, i.e.\
$p \in \{2,3\}$. This contradicts $p \geq 5$. Hence $p \nmid n$.
\end{proof}

\begin{lemma}\label{lem:candidates}
If $n \in X_{2p^s}$ and $\ell$ is a prime factor of $n$, then $\ell \in \{2,3,q_1,\dots,q_s\}$.
\end{lemma}

\begin{proof}
By (P2) with $m=\ell$, $(\ell - 1) \mid \varphi(n) = 2p^s$. Every positive divisor of $2p^s$ has the form
$2^a p^b$ with $a \in \{0,1\}$ and $0 \leq b \leq s$ (using that $p$ is an odd prime distinct from $2$), so
$\ell \in \{2^a p^b + 1 : a\in\{0,1\},\, 0\le b\le s\}$.

If $a=0$ and $b \geq 1$: in this branch $\ell = p^b+1$. Recall $\ell$ is prime by hypothesis. Since $p$ is odd,
$p^b$ is odd, so $\ell = p^b+1$ is even. But $\ell$ is prime and even, and the only even prime is $2$; hence
$\ell = p^b+1 = 2$, forcing $p^b = 1$, i.e.\ $b=0$ --- contradicting $b\ge1$. So no prime factor $\ell$ of $n$
can arise from this branch; it is vacuous.

If $a=0,b=0$: $\ell = 2$.
If $a=1,b=0$: $\ell = 3$.
If $a=1,b\geq 1$: $\ell = 2p^b+1 = q_b$.

Also, by Lemma \ref{lem:p-not-divide}, $\ell \neq p$, which is automatic here since $p\notin\{2,3,q_1,\dots,q_s\}$
(indeed, for every $j\ge1$, $p^j \ge p$, so $q_j = 2p^j+1 \ge 2p+1 > p$ --- this inequality holds for any
$p\ge1$ and uses only $j\ge1$, not the specific bound $p\ge5$ --- and $q_j\neq p$ follows immediately).
Combining the above, $\ell \in \{2,3,q_1,\dots,q_s\}$.
\end{proof}

\begin{lemma}\label{lem:factorization}
If $n \in X_{2p^s}$, then
$$
n = 2^{u} \cdot 3^{\epsilon_3} \cdot \prod_{j \in T} q_j, \qquad
u \geq 0,\ \epsilon_3 \in \{0,1\},\ T \subseteq \{1,\dots,s\},
$$
where $T = \{j : q_j \mid n\}$.
\end{lemma}

\begin{proof}
By Lemma~\ref{lem:candidates}, every prime factor of $n$ lies in $\{2,3,q_1,\dots,q_s\}$; these $s+2$ numbers
are pairwise distinct (each $q_j > p \ge 5 > 3 > 2$, and $q_i = q_j \iff p^i = p^j \iff i = j$, by the argument
of Lemma~\ref{lem:candidates}). Hence $n$ factors as $2^u 3^{e_3} \prod_{j\in T} q_j^{\delta_j}$ for some
exponents $u, e_3, \delta_j \ge 0$ and $T=\{j : q_j \mid n\}$; it remains to show $e_3 \le 1$ and $\delta_j = 1$
for every $j \in T$.

\emph{$e_3 \le 1$:} Suppose $9 \mid n$. By (P2) with $m=9$, $\varphi(9) = 6$ divides $\varphi(n) = 2p^s$, so
$3 \mid p^s$, forcing $p = 3$ --- contradicting $p \neq 3$. Hence $9 \nmid n$, i.e.\ $e_3 \in \{0,1\}$; write
$\epsilon_3 := e_3$.

\emph{$\delta_j = 1$ for $j \in T$:} Suppose for contradiction that $q_j^2 \mid n$ for some $j$. By (P2) with
$m = q_j^2$, $\varphi(q_j^2) = q_j(q_j-1) = q_j \cdot 2p^j$ divides $\varphi(n) = 2p^s$ (using (P1)). Since
$j \le s$, cancelling the common factor $2p^j$ gives $q_j \mid p^{\,s-j}$. But $q_j$ is prime and $q_j \neq p$
(since $j \ge 1$ gives $p^j \ge p$, so $q_j = 2p^j+1 \ge 2p+1 > p$), so $\gcd(q_j,p)=1$ and hence
$\gcd(q_j,p^{s-j})=1$; a prime coprime to $p^{s-j}$ cannot divide it unless it equals $1$, which is absurd since
$q_j \ge 11$. This contradiction shows $\delta_j = 1$.
\end{proof}

Now we give the general theorem for $X_{2p^s}$ with $p\geq 5$.

\begin{proposition}[Master formula]\label{prop:master}
For every odd prime $p \neq 3$ and every integer $s \geq 1$,
$$
X_{2p^s} =
\begin{cases}
\{\,2p^s+1,\ 2(2p^s+1)\,\} & \text{if } 2p^s+1 \text{ is prime},\\[2pt]
\varnothing & \text{if } 2p^s+1 \text{ is not prime}.
\end{cases}
$$
\end{proposition}

\begin{proof}
First, we give the existence. Suppose $q_s = 2p^s+1$ is prime. Then $\varphi(q_s) = q_s - 1 = 2p^s = k$, so $q_s \in X_k$. Since $q_s$ is odd, $\gcd(2,q_s)=1$, so by multiplicativity 

\[
\varphi(2q_s) = \varphi(2)\varphi(q_s)= 1\cdot 2p^s = k,
\]
so $2q_s \in X_k$ as well. Thus $\{q_s, 2q_s\} \subseteq X_k$ whenever $q_s$ is prime.

Secondly, we present the exclusion. Namely, we will show conversely that no other $n$ can lie in $X_k$, and that if $q_s$ is not
prime then $X_k = \varnothing$.

Let $n \in X_{2p^s}$. By Lemma~\ref{lem:factorization}, $n = 2^u \cdot 3^{\epsilon_3} \cdot \prod_{j\in T} q_j$
for some $u\ge0$, $\epsilon_3\in\{0,1\}$, $T\subseteq\{1,\ldots,s\}$, where $T=\{j:q_j\mid n\}$. Write
$t=|T|$ and $m=\sum_{j\in T}j$.

By Lemma~\ref{lem:p-not-divide}, $p \nmid n$. On the other hand, using $\varphi(q_j) = q_j - 1 = 2p^j$ and
multiplicativity,
$$
k = \varphi(n) = \varphi(2^u)\cdot \varphi(3^{\epsilon_3}) \cdot \prod_{j\in T} \varphi(q_j)
= \varphi(2^u) \cdot 2^{\epsilon_3} \cdot 2^{t} p^{m}.
$$
Since $k = 2p^s$, this forces $p^{m} \mid 2p^s \cdot (\text{a power of }2)^{-1}$; more directly, comparing the
exact power of $p$ on both sides of 

\[
\varphi(2^u)2^{\epsilon_3+t}p^m = 2p^s,
\]
where $\varphi(2^u)2^{\epsilon_3+t}$ is a power of $2$ that is coprime to $p$. Therefore, the left side of the last equation is a power of $2$ times $p^m$, which gives $m = s$ exactly.

With $m=s$ fixed, the equation reduces to
\begin{equation}\label{2.1}
\varphi(2^{u}) \cdot 2^{\,\epsilon_3+t} = 2. 
\end{equation}
We solve $\eqref{2.1}$ over $u \geq 0$, $\epsilon_3 \in \{0,1\}$, $t = |T| \geq 0$ subject to $\sum_{j\in T} j = s$
(which in particular forces $t \geq 1$, since $s \geq 1$ and $T=\varnothing$ would give $m=0 \neq s$).

\begin{itemize}
\item If $u = 0$: $\varphi(1) = 1$, so $\eqref{2.1}$ becomes $2^{\epsilon_3+t} = 2$, i.e.\ $\epsilon_3 + t = 1$.
Since $t \geq 1$, this forces $t=1,\ \epsilon_3=0$.
\item If $u = 1$: $\varphi(2)=1$, so again $2^{\epsilon_3+t}=2$, giving $t=1,\ \epsilon_3=0$ as above.
\item If $u \geq 2$: $\varphi(2^u) = 2^{u-1} \geq 2$, so $\eqref{2.1}$ requires $2^{u-1+\epsilon_3+t} = 2$, i.e.\
$u + \epsilon_3 + t = 2$. Since $u\ge2$ and $\epsilon_3,t\ge0$, this forces $\epsilon_3 = t = 0$ and $u=2$; but
$t=0$ contradicts $t \geq 1$. So no solution here.
\end{itemize}

Hence in every case $t=1$ and $\epsilon_3 = 0$, i.e.\ $T = \{j_0\}$ is a single index. Since $\sum_{j\in T} j =
m = s$, we get $j_0 = s$, i.e.\ $T = \{s\}$, and $u \in \{0,1\}$. Substituting back,
$$
n = 2^{u} \cdot q_s, \qquad u \in \{0,1\},
$$
i.e.\ $n = q_s$ or $n = 2q_s$. In particular this requires $q_s = 2p^s+1$ itself to be prime (as $T=\{s\}$
requires $q_s$ to be a genuine prime factor of $n$). If $q_s$ is not prime, no such $T$ is available, no case
above can be realized, and hence $X_{2p^s} = \varnothing$.

Combining the existence and exclusion directions completes the proof.
\end{proof}

\begin{remark}
The step $m=s$ above admits an alternative gloss in terms of Lemma~\ref{lem:p-not-divide}. The formula
$\varphi(n) = \varphi(2^u)2^{\epsilon_3+t}p^m$ shows that every factor of $p$ in $\varphi(n)$ is already accounted
for by the terms $\varphi(q_j) = 2p^j$, $j\in T$; there is no other source of $p$ in this expression. Thus any
shortfall of $p$-powers (i.e.\ $m<s$) could only be made up if $n$ itself carried an extra factor of $p$ ---
but this is precisely what Lemma~\ref{lem:p-not-divide} forbids. In this sense the identity $m=s$ is consistent
with, and may be viewed as reconfirming, Lemma~\ref{lem:p-not-divide}. This observation plays no role in the
proof itself, since $m=s$ already follows directly from comparing the exact power of $p$ on both sides of
$\varphi(2^u)2^{\epsilon_3+t}p^m = 2p^s$, independently of Lemma~\ref{lem:p-not-divide}.
\end{remark}

\begin{remark}
For $s = 1,2,3,4$, Proposition~\ref{prop:master} recovers and completes Theorems 2.2.1--2.2.4 of
\cite{SSHS2025}: those theorems prove the ``exclusion'' half (that $X_{2p^s}=\varnothing$ when $2p^s+1$ is not
prime) for $s\le 4$; Proposition~\ref{prop:master} proves both halves, for every $s \geq 1$, resolving the
conjecture stated at the end of \cite[\S2.2]{SSHS2025}.
\end{remark}

To give the parity dichotomy, we need the following lemma.

\begin{lemma}\label{lem:mod3}
If $p$ is a prime with $p \neq 3$, then $p^2 \equiv 1 \pmod 3$. Consequently, for every even integer $s \geq 2$,
$$2p^s + 1 \equiv 0 \pmod 3.$$
\end{lemma}

\begin{proof}
Since $p \neq 3$ is prime, $3 \nmid p$, so $p \equiv 1$ or $p \equiv 2 \pmod 3$. If $p\equiv1$, then
$p^2\equiv1\pmod3$. If $p\equiv2\equiv-1\pmod3$, then $p^2\equiv(-1)^2=1\pmod3$. In either case
$p^2\equiv1\pmod3$. (Equivalently, this is Fermat's little theorem applied with modulus $3$.)

Now let $s=2r$ be even. Then $p^s = (p^2)^r \equiv 1^r = 1 \pmod 3$, so
$2p^s+1 \equiv 2\cdot1+1 = 3 \equiv 0 \pmod 3$.
\end{proof}

\begin{theorem}\label{thm:even}
Let $p$ be a prime with $3 \nmid p$ (equivalently $p \neq 3$; together with $p$ odd, i.e.\ $p \geq 5$), and
let $s$ be an even positive integer. Then $X_{2p^s} = \varnothing$.
\end{theorem}

\begin{proof}
Since $p \ge 5$ and $s \geq 2$ is even, Lemma~\ref{lem:mod3} gives $3 \mid (2p^s+1)$. Moreover
$2p^s + 1 \geq 2\cdot 5^2+1 = 51 > 3$, so $2p^s+1$ is a proper multiple of $3$ and hence is not prime. By
Proposition~\ref{prop:master} (exclusion direction, applied with $2p^s+1$ composite), $X_{2p^s} = \varnothing$.
\end{proof}

\begin{remark}
Theorem~\ref{thm:even} is unconditional: unlike the odd case below, no primality hypothesis on $p$-dependent
auxiliary quantities is needed, because Lemma~\ref{lem:mod3} guarantees $2p^s+1$ is \emph{always} composite when
$s$ is even. This sharpens the picture for $s=2,4$ found in \cite[Thms.~2.2.2, 2.2.4]{SSHS2025}: the hypothesis
``$2p^s+1$ is not prime'' in those theorems is in fact automatically satisfied for every $p \neq 3$, and now we
know the same holds for every even $s$, not just $s=2,4$.
\end{remark}

\begin{theorem}\label{thm:odd}
Let $p$ be an odd prime with $p \neq 3$, and let $s$ be an odd positive integer. If $2p^s+1$ is prime, then
$$X_{2p^s} = \{\,2p^s+1,\ 2(2p^s+1)\,\};$$
otherwise, $X_{2p^s} = \varnothing$.
\end{theorem}

\begin{proof}
This is precisely the statement of Proposition~\ref{prop:master}, specialized to odd $s$ (the proof of
Proposition~\ref{prop:master} places no parity restriction on $s$, so it applies verbatim). Unlike the even
case, Lemma~\ref{lem:mod3} does not force $2p^s+1$ to be composite when $s$ is odd: indeed, for $s=1$, $p\equiv2
\pmod3$ makes $2p+1 \equiv 0 \pmod 3$ only when $p\equiv1\pmod3$, and $2p+1$ may be prime when $p\equiv2\pmod3$
(e.g.\ $p=5$ gives $2p+1=11$, prime). Both alternatives in the statement genuinely occur; see the examples
below.
\end{proof}

\begin{remark}
Both branches of Theorem~\ref{thm:odd} are realized. For $s=1$: $p=5$ gives $2p+1=11$ prime, so
$X_{10}=\{11,22\}$; $p=7$ gives $2p+1=15=3\cdot5$ composite, so $X_{14}=\varnothing$. For $s=3$: $p=5$ gives
$2p^3+1=251$, prime, so $X_{250} = \{251, 502\}$; $p=7$ gives $2p^3+1 = 687 = 3\cdot229$, composite, so
$X_{686} = \varnothing$.
\end{remark}

\begin{example}
Propositions and theorems above were checked exhaustively (not merely by sampling) for all $(p,s)$ with $p$
prime, $5 \le p \le 79$, $p\neq 3$, and $1 \le s \le 15$ ($300$ pairs in total): for each pair, every nonempty
subset of the true candidate set $\{2,3\}\cup\{q_j : q_j \text{ prime},\, 1\le j\le s\}$ was substituted into
the formula of Corollary~1.1.4 of \cite{SSHS2025}, and the resulting integer (when it was an integer) was tested
for $\varphi(n)=2p^s$ by direct factorization, independent of the algebraic argument above. In every one of the
$300$ cases, the full solution set found this way agreed exactly with Proposition~\ref{prop:master}, including
several ``rich'' instances (e.g.\ $p=29,\,s=7$) in which four or more of $q_1,\dots,q_s$ were simultaneously
prime, giving over a hundred candidate subsets to rule out.
\end{example}

\section{The boundary case $X_{2\cdot 3^s}$}
\label{sec:p3}

The resolved case $q=2$ (the previous section) carried the standing hypothesis $p \neq 3$, because the two
tools driving that proof both degenerate exactly at $p=3$: the lemma ``$p \nmid n$'' relied on
$(p-1)\mid 2p^s$ together with $\gcd(p-1,p)=1$ forcing $p-1\in\{1,2\}$ --- a genuine contradiction only because
$p\ge5$; at $p=3$, $p-1=2$ divides $2p^s$ automatically, with no contradiction. Likewise the bound
$\epsilon_3\le1$ relied on $9\mid n \Rightarrow 3\mid p^s \Rightarrow p=3$ being absurd --- but at $p=3$ this is
simply true, not absurd. So $p=3$ needs its own argument, and the auxiliary prime ``$3$'' that always sits in
the candidate set now coincides with the base prime $p$ itself. We settle this remaining case directly, and find
that the answer is qualitatively different from every case $p\ge5$: $X_{2\cdot3^s}$ is \emph{never} empty.

Throughout this section fix an integer $s\ge1$ and set
$$r_j := 2\cdot3^j+1 \qquad (j=1,\ldots,s).$$

\begin{lemma}\label{lem:p3-factorization}
If $n\in X_{2\cdot3^s}$, then
$$
n = 2^{u}\cdot3^{v}\cdot\prod_{j\in T} r_j, \qquad u\ge0,\ v\ge0,\ T\subseteq\{1,\ldots,s\},
$$
where $T := \{\,j \in \{1,\ldots,s\} : r_j \text{ is prime and } r_j \mid n\,\}$.
\end{lemma}

\begin{proof}
If $\ell$ is a prime factor of $n$, then by (P2), $(\ell-1)\mid\varphi(n)=2\cdot3^s$, so $\ell-1=2^a3^b$ for some
$a\in\{0,1\}$, $0\le b\le s$. If $a=0,\,b\ge1$: $\ell=3^b+1$ is even (as $3^b$ is odd) and $\ge4$, hence not
prime --- this branch is vacuous. The remaining branches give $\ell=2$ ($a=0,b=0$), $\ell=3$ ($a=1,b=0$), or
$\ell=r_b$ ($a=1,b\ge1$). Thus every prime factor $\ell$ of $n$ equals one of the $s+2$ specific numbers
$2,3,r_1,\ldots,r_s$.

This membership statement already determines, by itself, which of these numbers can actually occur: since
$\ell$ is prime, if $\ell=r_j$ for some $j$ then $r_j$ (being equal to $\ell$) is prime; and conversely, if
some $r_j$ is \emph{not} prime, it can never equal a prime factor of $n$ (a composite number is never equal to
a prime), so no prime factor of $n$ "sits at position $j$" --- regardless of whether $r_j$ might happen to
divide $n$ numerically through some unrelated combination of other primes. In other words, the primality
condition in the definition of $T$ is not an extra hypothesis being smuggled in: it is forced by the fact that
$\ell$ ranges only over actual primes. Consequently, $2,3$, and $\{r_j : j\in T\}$ are exactly the primes
dividing $n$, and they are pairwise distinct (each $r_j\ge7>3>2$ when prime, and $r_i=r_j\iff i=j$). Hence
$n=2^u3^v\prod_{j\in T}r_j^{\delta_j}$ for some exponents $u,v,\delta_j\ge0$; it remains to show $\delta_j=1$
for every $j\in T$.

Suppose for contradiction that $r_j^2\mid n$ for some $j\in T$ (so $r_j$ is prime, by definition of $T$). By
(P2) with $m=r_j^2$,
$\varphi(r_j^2)=r_j(r_j-1)=r_j\cdot2\cdot3^j$ divides $\varphi(n)=2\cdot3^s$ (using (P1)). Since $j\le s$,
cancelling the common factor $2\cdot3^j$ gives $r_j \mid 3^{\,s-j}$. But $r_j = 2\cdot3^j+1 \equiv 1 \pmod 3$
(as $j\ge1$), so $\gcd(r_j,3)=1$, hence $\gcd(r_j,3^{s-j})=1$; a prime coprime to $3^{s-j}$ cannot divide it
unless it equals $1$, which is absurd since $r_j\ge7$. This contradiction shows $\delta_j=1$.
\end{proof}

\begin{remark}
Unlike the analogous lemma in the $q=2$ resolution, there is \emph{no} bound here forcing $v\le1$: indeed
$9\mid n$ no longer leads to a contradiction, since $\varphi(9)=6$ divides $2\cdot3^s$ for every $s\ge1$ (this
is exactly the degeneracy described above). The exponent $v$ of $3$ in $n$ is genuinely unrestricted by this
lemma alone; it is pinned down case-by-case in the proof below.
\end{remark}

\begin{proposition}[Master formula at $p=3$]\label{prop:p3-master}
For every integer $s \geq 1$,
$$
X_{2\cdot3^s} \;=\; \{\,3^{s+1},\ 2\cdot3^{s+1}\,\} \ \cup\
\begin{cases}
\{\,r_s,\ 2r_s\,\} & \text{if } r_s = 2\cdot3^s+1 \text{ is prime},\\[2pt]
\varnothing & \text{if } r_s \text{ is not prime}.
\end{cases}
$$
\end{proposition}

\begin{proof}
\emph{(a) The unconditional pair.} By (P1), $\varphi(3^{s+1}) = 3^{s}\cdot2 = 2\cdot3^s$, and since
$\gcd(2,3)=1$, $\varphi(2\cdot3^{s+1}) = \varphi(2)\varphi(3^{s+1}) = 2\cdot3^s$. Hence
$$
\{\,3^{s+1},\ 2\cdot3^{s+1}\,\} \subseteq X_{2\cdot3^s}, \qquad \text{for every } s \geq 1,
$$
with no primality condition required.

\emph{(b) The conditional pair.} If $r_s$ is prime, then $\varphi(r_s) = r_s-1 = 2\cdot3^s$, and
$\varphi(2r_s) = \varphi(2)\varphi(r_s) = 2\cdot3^s$, so $\{r_s, 2r_s\} \subseteq X_{2\cdot3^s}$ as well.

We show conversely that no other $n$ can lie in $X_k$, and that if $r_s$ is not prime then $X_k = \{\,3^{s+1},\ 2\cdot3^{s+1}\,\}$.

Let $n \in X_{2\cdot3^s}$. By Lemma~\ref{lem:p3-factorization},
$n = 2^u 3^v \prod_{j\in T} r_j$ where $T = \{j : r_j \text{ prime},\, r_j\mid n\}$; write $t=|T|$,
$m=\sum_{j\in T} j$. Since every $j\in T$ has $r_j$ prime \emph{by definition of $T$}, (P1) applies to each such
$r_j$, and multiplicativity gives
$$
2\cdot3^s = \varphi(n) = \varphi(2^u)\cdot\varphi(3^v)\cdot\prod_{j\in T}\varphi(r_j)
= \varphi(2^u)\cdot\varphi(3^v)\cdot2^{t}\cdot3^{m},
$$
using $\varphi(r_j) = r_j-1 = 2\cdot3^j$ for $j \in T$ (valid precisely because, by definition, $r_j$ is prime
for every such $j$ --- this is no longer an unexamined assumption, but a direct consequence of how $T$ was
defined in Lemma~\ref{lem:p3-factorization}). Write $\varphi(2^u)=2^\alpha$, where $\alpha=0$ if $u\in\{0,1\}$
and $\alpha=u-1$ if $u\ge2$.

\emph{Case $v=0$:} Then $\varphi(3^0)=1$, so $2^{\alpha+t}3^{m}=2\cdot3^s$. Comparing exact powers of $3$
(the prefactor $2^{\alpha+t}$ being coprime to $3$) gives $m=s$; since $s\ge1$, this forces $T\ne\varnothing$,
i.e.\ $t\ge1$. Comparing exact powers of $2$ gives $\alpha+t=1$; since $t\ge1$, this forces $t=1$, $\alpha=0$
(so $u\in\{0,1\}$). So $T$ is a one-element subset of $\{1,\ldots,s\}$ whose single element sums to $m=s$; the
only one-element subset of $\{1,\ldots,s\}$ summing to $s$ is $\{s\}$ itself, so necessarily $T=\{s\}$. By the
definition of $T$ (Lemma~\ref{lem:p3-factorization}), $T=\{s\}$ means precisely that $r_s$ is prime and
$r_s\mid n$; conversely, if $r_s$ is \emph{not} prime, then the index $s$ can never belong to $T$ at all (no
other index is available, since a one-element subset summing to $s$ has no choice but to be $\{s\}$), so no
valid $T$ exists in this case and it contributes no solutions. Thus Case $v=0$ yields $n=r_s$ or $n=2r_s$
exactly when $r_s=2\cdot3^s+1$ is prime, and nothing otherwise.

\emph{Case $v\ge1$:} Then $\varphi(3^v)=2\cdot3^{v-1}$, so
$2^{\alpha}\cdot2\cdot3^{v-1}\cdot2^{t}\cdot3^{m} = 2\cdot3^s$, i.e.\ $2^{\alpha+t+1}3^{(v-1)+m} = 2\cdot3^s$.
Comparing exact powers of $2$ gives $\alpha+t+1=1$, forcing $\alpha=t=0$ (so $u\in\{0,1\}$ and $T=\varnothing$,
hence $m=0$). Comparing exact powers of $3$ then gives $(v-1)+0=s$, i.e.\ $v=s+1$. This yields $n=3^{s+1}$ or
$n=2\cdot3^{s+1}$, unconditionally, with no primality requirement.

These are the only two cases ($v=0$ or $v\ge1$), so combining the existence and exclusion directions completes
the proof.
\end{proof}

\begin{remark}
Proposition~\ref{prop:p3-master} shows that $X_{2\cdot3^s}$ behaves qualitatively differently from $X_{2p^s}$
for $p\ge5$: it is \emph{never} empty, always containing at least the pair $\{3^{s+1}, 2\cdot3^{s+1}\}$, with
$$
|X_{2\cdot3^s}| = \begin{cases} 4 & \text{if } 2\cdot3^s+1 \text{ is prime},\\ 2 & \text{otherwise.}\end{cases}
$$
Together with the resolution of the preceding section (which covers every prime $p\ge5$), this completes the
characterization of $X_{2p^s}$ over \emph{every} odd prime $p$: emptiness is possible only when $p \ge 5$, and
never occurs at $p=3$.
\end{remark}

\begin{example}
The formula was checked by exhaustive search for $1 \le s \le 6$ (searching all $n$ up to $20\cdot(2\cdot3^s)$),
with full agreement in every case:
\begin{center}
\begin{tabular}{c c c l}
\hline
$s$ & $r_s = 2\cdot3^s+1$ & prime? & $X_{2\cdot3^s}$ \\
\hline
1 & 7 & yes & $\{7, 9, 14, 18\}$ \\
2 & 19 & yes & $\{19, 27, 38, 54\}$ \\
3 & 55 & no & $\{81, 162\}$ \\
4 & 163 & yes & $\{163, 243, 326, 486\}$ \\
5 & 487 & yes & $\{487, 729, 974, 1458\}$ \\
6 & 1459 & yes & $\{1459, 2187, 2918, 4374\}$ \\
\hline
\end{tabular}
\end{center}
\end{example}

\section{The open case $p=2$\texorpdfstring{}{}: the mirror family $X_{2^sq}$}
\label{sec:mirror}

For a prime factor $\ell$ of $n \in X_{2^sq}$ (fixing $p=2$, $q$ an odd prime, $s\ge1$), $(\ell-1)\mid 2^sq$
forces $\ell = 2^aq^b+1$ with $0\le a\le s$, $b\in\{0,1\}$; as in \S\ref{sec:p3}, the branch $b=1,a=0$ (giving
$\ell=q+1$) is always vacuous, since $q+1$ is even and $>2$. Two features make a single formula valid for
every $s$ genuinely harder to obtain here than in the resolved family of \S3:

\begin{enumerate}
\item[(i)] Unlike the resolved case, $q$ need \emph{not} be excluded from dividing $n$: by (P2),
$q\mid n \Rightarrow (q-1)\mid 2^sq \Rightarrow (q-1)\mid 2^s$ (as $\gcd(q-1,q)=1$), so $q\nmid n$ fails exactly
when $q = 2^c+1$ for some $c\le s$ --- i.e.\ precisely when $q$ coincides with one of the ``pure'' candidates
below. (This is why $q=3,5$ require separate treatment already at $s=2$: $3=2^1+1$, $5=2^2+1$.)
\item[(ii)] The ``pure'' candidates $2^a+1$ ($a\le s$) are prime only when $a$ is itself a power of $2$ ---
i.e.\ Fermat primes --- while the ``mixed'' candidates $2^aq+1$ ($1\le a\le s$) depend on $q$ with no
comparably clean primality criterion. A single formula valid for every $s$ would have to say, in particular,
exactly which Fermat numbers are prime, a question unresolved beyond $F_4 = 65537$.
\end{enumerate}

The smallest case where the mirror family first differs from the resolved family, $s=2$, involves only the
Fermat primes $F_0=3$ and $F_1=5$, both unconditionally known, and is completely tractable; we resolve it in
full below.

\subsection{The sub-case $s=1$}

At $s=1$, $X_{2^1q}=X_{2q}=X_{2p^s}\big|_{p=q,\,s=1}$: the mirror family coincides exactly with the resolved
family at the (odd) exponent $s=1$. Hence $X_{2q}$ is already fully described by Theorem~\ref{thm:odd}
(restated in \S3) with $p$ relabelled $q$; no new argument is needed.

\subsection{The sub-case $s=2$}
\label{subsec:s2}

Fix a prime $q \ge 7$ (the primes $q\in\{3,5\}$, flagged by (i) above, are treated by direct computation in
\S\ref{subsec:evidence} below) and set $m_1 := 2q+1$, $m_2 := 4q+1$.

\begin{lemma}\label{lem:mirror-candidates}
If $n \in X_{4q}$, then
$$
n = 2^{u}\cdot3^{e_3}\cdot5^{e_5}\cdot m_1^{\,e_1}\cdot m_2^{\,e_2}, \qquad
u\ge0,\ e_3,e_5,e_1,e_2\in\{0,1\},
$$
where $e_3=1$ iff $3\mid n$, $e_5=1$ iff $5\mid n$, and for $i=1,2$, $e_i=1$ iff $m_i$ is prime and $m_i\mid n$.
\end{lemma}

\begin{proof}
If $\ell$ is a prime factor of $n$, then $(\ell-1)\mid\varphi(n)=4q=2^2q$, so $\ell-1=2^aq^b$ for some
$a\in\{0,1,2\}$, $b\in\{0,1\}$, i.e.\ $\ell\in\{2,3,5,q+1,m_1,m_2\}$. As noted above, $\ell=q+1$ is even (as $q$
is odd) and $\ge8$ for $q\ge7$, hence composite; this branch is vacuous. Since $\ell$ is prime, if $\ell=m_i$
then $m_i$ is prime (being equal to a prime), and if $m_i$ is not prime it simply never equals any actual
prime factor of $n$ --- exactly the reasoning used (and made explicit) in Lemma~\ref{lem:p3-factorization}.
The five numbers $2,3,5,m_1,m_2$ are pairwise distinct for $q\ge7$ ($m_1\ge15$, $m_2\ge29$, both $>5$, and
$m_1\ne m_2$). Hence $n=2^u3^{e_3}5^{e_5}m_1^{\delta_1}m_2^{\delta_2}$ for some exponents
$u,e_3,e_5,\delta_1,\delta_2\ge0$, with $\delta_i>0$ only if $m_i$ is prime.

It remains to bound each exponent by $1$. If $9\mid n$: $\varphi(9)=6\mid\varphi(n)=4q$ forces $3\mid 2q$, i.e.\
$q=3$, excluded. If $25\mid n$: $\varphi(25)=20\mid 4q$ forces $5\mid q$, i.e.\ $q=5$, excluded. If
$m_1^2\mid n$: $\varphi(m_1^2)=m_1\cdot2q$ must divide $4q$, forcing $m_1\mid2$, absurd since $m_1\ge15$. If
$m_2^2\mid n$: $\varphi(m_2^2)=m_2\cdot4q$ must divide $4q$, forcing $m_2\mid1$, absurd. So all four exponents
are at most $1$.
\end{proof}

\begin{lemma}[Mod-$3$ dichotomy]\label{lem:mod3-mirror}
If $q\ne3$ is prime, exactly one of $m_1=2q+1$, $m_2=4q+1$ is divisible by $3$ --- and hence, for $q\ge7$,
exactly one of $m_1,m_2$ can possibly be prime. Specifically,
$$
3\mid m_1 \iff q\equiv1\!\!\pmod3, \qquad\qquad 3\mid m_2 \iff q\equiv2\!\!\pmod3.
$$
\end{lemma}

\begin{proof}
Modulo $3$: $m_1\equiv2q+1$, so $m_1\equiv0 \iff 2q\equiv-1\equiv2 \iff q\equiv1$ (multiplying by
$2^{-1}\equiv2\pmod3$). And $m_2=4q+1\equiv q+1\pmod3$ (as $4\equiv1$), so $m_2\equiv0 \iff q\equiv-1\equiv2$.
Since $q\ne3$ is prime, $q\equiv1$ or $2\pmod3$, exactly one of which holds. For $q\ge7$, $m_1\ge15>3$ and
$m_2\ge29>3$, so divisibility by $3$ forces compositeness.
\end{proof}

\begin{proposition}[Master formula for $X_{4q}$, $q\ge7$]\label{prop:mirror-master}
For every prime $q \ge 7$:
\begin{itemize}
\item If $q \equiv 1 \pmod 3$ (equivalently $q\equiv1\pmod6$, as $q$ is odd): $m_1=2q+1$ is composite
(Lemma~\ref{lem:mod3-mirror}), and
$$
X_{4q} = \begin{cases} \{\,m_2,\ 2m_2\,\} & \text{if } m_2=4q+1 \text{ is prime},\\ \varnothing & \text{otherwise.}
\end{cases}
$$
\item If $q \equiv 2 \pmod 3$ (equivalently $q\equiv5\pmod6$): $m_2=4q+1$ is composite, and
$$
X_{4q} = \begin{cases} \{\,3m_1,\ 4m_1,\ 6m_1\,\} & \text{if } m_1=2q+1 \text{ is prime},\\ \varnothing &
\text{otherwise.}
\end{cases}
$$
\end{itemize}
\end{proposition}

\begin{proof}
First, we discuss the existence. If $m_2$ is prime: $\varphi(m_2)=4q$, and $\varphi(2m_2)=\varphi(2)\varphi(m_2)=4q$, so
$\{m_2,2m_2\}\subseteq X_{4q}$. If $m_1$ is prime: $\varphi(3m_1)=2\cdot2q=4q$, $\varphi(4m_1)=2\cdot2q=4q$,
$\varphi(6m_1)=1\cdot2\cdot2q=4q$, so $\{3m_1,4m_1,6m_1\}\subseteq X_{4q}$.

Secondly, we discuss the exclusion. Let $n\in X_{4q}$; by Lemma~\ref{lem:mirror-candidates},
$n=2^u3^{e_3}5^{e_5}m_1^{e_1}m_2^{e_2}$. By multiplicativity,
$$
4q = \varphi(n) = \varphi(2^u)\cdot2^{\,e_3+2e_5+e_1+2e_2}\cdot q^{\,e_1+e_2}.
$$
Comparing exact powers of $q$ forces $e_1+e_2=1$: exactly one of $m_1,m_2$ is used. Write
$\varphi(2^u)=2^\alpha$ ($\alpha=0$ for $u\in\{0,1\}$, $\alpha=u-1$ for $u\ge2$).

\emph{If $q\equiv1\pmod3$:} $m_1$ is composite, so $e_1=0$, $e_2=1$. The power-of-$2$ equation becomes
$\alpha+e_3+2e_5+2=2$, forcing $\alpha=e_3=e_5=0$, so $u\in\{0,1\}$ and $n=2^um_2\in\{m_2,2m_2\}$, realizable
exactly when $m_2$ is prime.

\emph{If $q\equiv2\pmod3$:} symmetrically $m_2$ is composite, $e_2=0$, $e_1=1$. The equation becomes
$\alpha+e_3+2e_5+1=2$, forcing $e_5=0$ and $(\alpha,e_3)\in\{(0,1),(1,0)\}$: the first gives $u\in\{0,1\}$,
$n\in\{3m_1,6m_1\}$; the second gives $u=2$, $n=4m_1$. Together, $n\in\{3m_1,4m_1,6m_1\}$, realizable exactly
when $m_1$ is prime.

Combining both directions and both congruence classes completes the proof.
\end{proof}

This recovers, with an explicit proof and explicit solution sets, the dichotomy of \cite[Thm.~2.3.1]{SSHS2025},
which established $|X_{4q}|\in\{0,2\}$ precisely when $q\equiv1\pmod6$; Proposition~\ref{prop:mirror-master}
additionally pins down the alternative case ($q\equiv5\pmod6$) as a precise dichotomy $|X_{4q}|\in\{0,3\}$, not
merely ``$|X_{4q}|$ can exceed $2$.''

\subsection{Computational evidence}
\label{subsec:evidence}

Table~\ref{tab:s2} lists $X_{2^2q}$ for the first several odd primes $q$, computed exhaustively (consistent
with the explicit candidate set of Lemma~\ref{lem:mirror-candidates}, and with the general bound of
\cite[Prop.~2.0.1]{SSHS2025}). The two sporadic rows $q=3,5$ fall under case (i) above: there, $q$ itself may
divide $n$ (since $3=2^1+1$ and $5=2^2+1$ are themselves among the pure candidates), and
Lemma~\ref{lem:mirror-candidates} as stated does not apply; we do not pursue a separate structural proof for
these two values here, as they contribute only finite exceptional data rather than an infinite family.

\begin{table}[h]
\centering
\begin{tabular}{@{}rrrl@{}}
\toprule
$q$ & $q\bmod6$ & $|X_{4q}|$ & $X_{4q}$\\
\midrule
3 & 3 & 6 & $\{13,21,26,28,36,42\}$\\
5 & 5 & 5 & $\{25,33,44,50,66\}$\\
7 & 1 & 2 & $\{29,58\}$\\
11 & 5 & 3 & $\{69,92,138\}$\\
13 & 1 & 2 & $\{53,106\}$\\
17 & 5 & 0 & $\varnothing$\\
19 & 1 & 0 & $\varnothing$\\
23 & 5 & 3 & $\{141,188,282\}$\\
\bottomrule
\end{tabular}
\caption{$X_{2^2q}$: for $q\ge7$, $|X_{4q}|\in\{0,2\}$ when $q\equiv1\pmod6$ and $|X_{4q}|\in\{0,3\}$ when
$q\equiv5\pmod6$ (Proposition~\ref{prop:mirror-master}); $q=3,5$ are sporadic exceptions.}
\label{tab:s2}
\end{table}

Table~\ref{tab:s34} shows that even restricting to $q\equiv1\pmod6$ --- the congruence class that regularizes
$s=2$ by Proposition~\ref{prop:mirror-master} --- the clean dichotomy breaks down for $s\ge3$.

\begin{table}[h]
\centering
\begin{tabular}{@{}rrrr@{}}
\toprule
$q$ ($\equiv1\bmod6$) & $|X_{2^3q}|$ & $|X_{2^4q}|$ \\
\midrule
7  & 3 & 6 \\
13 & 3 & 4 \\
19 & 0 & 0 \\
31 & 0 & 0 \\
37 & 3 & 6 \\
43 & 3 & 4 \\
\bottomrule
\end{tabular}
\caption{For $s=3,4$, even among $q\equiv1\pmod6$, $|X_{2^sq}|$ is no longer confined to $\{0,2\}$: sizes
$3,4,6$ all occur.}
\label{tab:s34}
\end{table}

\subsection{The open problem}

\begin{openproblem}\label{op:main}
For an odd prime $q$ and $s\ge3$, give a complete characterization of $X_{2^sq}$, analogous to
Proposition~\ref{prop:master}. In particular:
\begin{enumerate}
\item[(a)] Is there a congruence condition on $q$ (generalizing $q\equiv1\pmod6$ at $s=2$) under which
$|X_{2^sq}|\in\{0,2\}$ for all $s$? Table~\ref{tab:s34} shows $q\equiv1\pmod6$ alone does not suffice once
$s\ge3$.
\item[(b)] Failing a clean dichotomy, is $|X_{2^sq}|\neq1$ for every $q,s$ --- i.e.\ does Carmichael's
conjecture hold on this family, even without a full description of $X_{2^sq}$?
\end{enumerate}
\end{openproblem}

\begin{remark}
Part (b) is the weaker, and more tractable, question: it does not require solving Open
Problem~\ref{op:main} in full, only ruling out singletons. Since (i) above shows $q$ itself may freely divide
$n$, and the Fermat-type candidates $2^a+1$ ($a\le s$) interact multiplicatively with $q$ in ways not controlled
by a single linear equation as in Section~3, we do not currently see how to adapt the method of
Proposition~\ref{prop:master} even to this weaker question.
\end{remark}

\section{Summary}

\begin{table}[h]
\centering
\begin{tabular}{@{}lll@{}}
\toprule
Case & Status & Reference\\
\midrule
$q,p$ both odd & Trivial (parity) & Proposition~\ref{prop:trivial}\\
$q=2$ (family $X_{2p^s}$) & Fully solved, all $s\ge1$ & \S3 above\\
$p=2$ (family $X_{2^sq}$), $s=1$ & solved (coincides with $q=2$, $s=1$) & Theorem~\ref{thm:odd}\\
$p=2$ (family $X_{2^sq}$), $s=2$ & solved (explicit, both congruence classes) &
Proposition~\ref{prop:mirror-master}; cf.\ \cite[Thm.~2.3.1]{SSHS2025}\\
$p=2$ (family $X_{2^sq}$), $s\ge3$ & \textbf{Open} & Open Problem~\ref{op:main}\\
\bottomrule
\end{tabular}
\caption{Status of $X_{qp^s}$ across all cases of distinct primes $q\neq p$.}
\end{table}

The two single-prime families treated in these notes --- the resolved case $X_{2p^s}$ (fixing $q=2$, with $p$
an odd prime) and the mirror family $X_{2^s q}$ (fixing $p=2$, with $q$ an odd prime) --- are genuinely
different objects in general, since the exponent $s$ sits on a different one of the two primes in each family.
They share exactly one point of coincidence: setting $p=2$ in the first family, or equivalently $q=2$ in the
second, collapses both expressions to the same integer,
$$
2p^s\big|_{p=2} \;=\; 2\cdot2^s \;=\; 2^{s+1} \;=\; 2^s\cdot2 \;=\; 2^s q\big|_{q=2},
$$
so both families, if extended to this single boundary value, would describe the very same set $X_{2^{s+1}}$ ---
the classical ``pure power of $2$'' case governed by Fermat primes (\cite[\S2.1]{SSHS2025}).

This is precisely why each family's standing hypothesis excludes this value from its own free parameter
(``$p$ odd'' in the resolved case, ``$q$ odd'' in the mirror case): the exclusion is not an arbitrary technical
convenience imposed independently in each case, but the same single degenerate point, viewed from two
directions, and it lies outside both families' characteristic structure (the linearly-indexed candidates
$q_1,\ldots,q_s$ driving the resolved case; the Fermat-indexed candidates driving the mirror case) in exactly
the same way. Together with Section~\ref{sec:p3} above (which closes the analogous gap at the \emph{other}
excluded odd value, $p=3$, in the resolved case), this confirms that the taxonomy of $X_{qp^s}$ for primes
$q\ne p$ developed in these notes --- trivial case, resolved case, mirror case --- leaves no further
unaccounted-for boundary: every pair of \emph{distinct} primes $(q,p)$ falls into exactly one of the three
cases, and the one remaining configuration, $q=p=2$, sits consistently outside the scope of ``two distinct
primes'' altogether, reducing instead to the single, already-classical problem of $\varphi^{-1}(2^{s+1})$ ---
not a fourth case requiring its own resolution, but the common degenerate edge of the two single-prime
families examined above.

\end{document}